\RequirePackage{fix-cm}
\documentclass{elsarticle}
\usepackage{amsmath}
\usepackage{amssymb}
\usepackage{lineno,hyperref}
\modulolinenumbers[5]
\usepackage[utf8]{inputenc}

\journal{}

\usepackage{amsmath}
\begin{document}

\begin{frontmatter}

\title{Multiparametric robust solutions for combinatorial problems  with parameterized locally budgeted uncertainty} 

\author{Alejandro Crema\\
 Escuela de Computaci\'on\\ Facultad de Ciencias\\ Universidad Central de Venezuela}
\newtheorem{teorema}{Teorema}
\newtheorem{lemma}{Lemma}
\newtheorem{corolario}{Corolario}
\newtheorem{ejemplo}{Ejemplo}
\newtheorem{capitulo}{Chapter}
\newtheorem{proposicion}{Proposition}
\newtheorem{remark}{Remarks}
\bibliographystyle{elsarticle-num}


\begin{abstract}
In this paper we studied combinatorial problems with parameterized locally budgeted uncertainty. We are looking for a solutions set such that for any parameters vector there exists a solution in the set with robustness near optimal. The algorithm consists of applying a multiparametric algorithm to obtain a near optimal multiparametric solution relative to the objective function for a combinatorial problem defined to find a robust solution for parameters fixed. As far as we know this is the first algorithm presented to do that task. Computational experience is presented to shortest path  and $p$-medians problems.\\

\end{abstract}

\end{frontmatter}

\section{Introduction}

Let $X \subseteq {\lbrace 0,1 \rbrace}^n$ with $X \neq \emptyset$. Let us suppose that $X$ is 0-1-Mixed Integer Linear Programming (0-1-MILP) representable and let $P(c)$ be a combinatorial problem in $x$, parameterized in $c \in \mathbb{R}^n_+$, defined as follows:
\begin{flalign*}
\underset{x \in X}{\min}&\;\;c^t x&P(c)&\
\end{flalign*}

Data uncertainty appears in many optimization problems. There are several options to model the uncertainty of the cost vector of a combinatorial optimization problem. Many examples may be seen in \cite{ChasGoe2016},\cite{ChasGoe2018},\cite{BuchKurtz2018},\cite{BenHerVial2015}.\\

Let $\Omega \subseteq \mathbb{R}^K_+$ and let us suppose that $\Omega$ is 0-1-MILP-representable, let $\mathcal{U} \in \mathbb{R}^K_+$ such that $\Gamma \leq \mathcal{U}$ for all $\Gamma \in \Omega$,  let $\underline{c},d \in \mathbb{R}^n_+$, let $[q]=\lbrace 1, \cdots,q\rbrace$ for all $q \in \mathbb{N}$ and let $\lbrace P_1,\cdots,P_K\rbrace$ be a partition of $[n]$. In this paper we consider locally budgeted uncertainty sets parameterized with $\Gamma \in \Omega$ as follows (\cite{GoeLen2021}):
\[\Lambda(\Gamma) = \lbrace c \in \mathbb{R}^n_+: c = \underline{c} + \lambda,\;\underset{j \in P_k}{\sum} \lambda_j \leq \Gamma_k\;\forall k \in [K],\;\lambda_j \in [0,d_j]\;\forall j \in [n]   \rbrace\]

Some examples of practical interest for $\Omega$ are the interval case, the line segment case and the budgeted case as follows:

\begin{itemize}
	\item Interval case: let $\mathcal{L} \in \mathbb{R}^K_+$ with $\mathcal{L} \leq \mathcal{U}$ and let $\Omega = [\mathcal{L},\mathcal{U}] = \lbrace \Gamma \in \mathbb{R}^K_+: \mathcal{L}  \leq \Gamma \leq \mathcal{U} \rbrace$
	\item Line segment case: let $\Gamma^0 \in \mathbb{R}^K_+$, let $[\underline{\alpha},\overline{\alpha}] \subseteq [0,1]$ and let $\Omega = \lbrace \Gamma \in \mathbb{R}^K_+:\Gamma = \alpha  \Gamma^0,\;\alpha \in [\underline{\alpha},\overline{\alpha}] \rbrace$ with $\mathcal{U} = \overline{\alpha}  \Gamma^0$
	\item Budgeted case: let $\underline{\Gamma},\mathcal{D} \in \mathbb{R}^K_+$, let $\Delta \in \mathbb{R}_+$ and let $\Omega =\lbrace \Gamma \in \mathbb{R}^K_+:\Gamma = \underline{\Gamma} + \beta,\;\;\underset{k \in [K]}{\sum} \beta_k \leq \Delta,\;\beta_k \in [0,\mathcal{D}_k]\;\forall k \in [K] \rbrace$ with  $\mathcal{U} = \underline{\Gamma} + \mathcal{D}$  
	
\end{itemize}

In recent decades Robust (\cite{ChasGoe2016}, \cite{BuchKurtz2018}), Stochastic (\cite{LiLiu2016}), Multiparametric \cite{OveDianNasPaSunAvPis2016}) and Fuzzy programming (\cite{CarlFu2002}) approaches have been developed to deal with such uncertainties. In this paper we use the Robust approach.\\

There are several robust optimization concepts you may select. Some examples are (\cite{ChasGoe2016}): classic robustness, absolute or relative regret robustness, adjustable robustness, recoverable robustness, light robustness, soft robustness, lexicographic $\alpha$-robustness, recovery-to-optimality, or similarity-based robustness. In this paper we consider the classic robustness.\\

Let $\Gamma \in \Omega $ and let $x \in X$. Let $W(x,\Gamma)$ be a Linear Programming (LP) problem in $c$ defined to compute the {\it robustness of $x$ for $\Gamma$}, as follows: 
\begin{flalign*}  
\underset{c \in \Lambda(\Gamma)}  {\max}&\;\;c^t x&W(x,\Gamma)&\
\end{flalign*}

In the rest of the paper if $S$ is an optimization problem then $v(S)$ is its optimal value.\\ 

Classical robust approach is to find $x \in X$ with minimal robustness for $\Gamma$ by solving the following problem in $x$:
\begin{flalign*}
\underset{x \in X} {\min}&\;\;v(W(x,\Gamma)) = \underset{x \in X} {\min}\;\left ( \underset{c\in \Lambda(\Gamma)}{\max}\;c^t x \right )&R(\Gamma)\
\end{flalign*}

If $x_R(\Gamma)$ is an optimal solution for $R(\Gamma)$ then $x_R(\Gamma)$ is a {\it robust solution for $\Gamma$}.\\

In this paper we consider the question how robust solutions change when the uncertainty set changes. In our study the structure of the uncertainty is always the same but some parameters change. Our goal is to compute a set  ${\lbrace x^i \rbrace}_1^r \subseteq X$ such
that for any $\Gamma \in \Omega$, there exists $i(\Gamma) \in [r]$ such that $x^{i(\Gamma)}$ is near optimal for $R(\Gamma)$. A pionner work with a single parameter $\lambda$ controlling the size of the uncertainty set, defined as $\Lambda(\lambda) = \underline{c} + \lambda B$ with $B$ a convex set containing the origin, may be seen in \cite{ChasGoe2018}.\\

Beyond the theoretical interest of a work like the one presented, the practical idea behind this approach is to find a set of solutions and choose the best one each time a new scenario $c$ appears. As an application, imagine a public service system designed based on the shortest path (or $p$-medians) problem. Each time that changes the current situation a new path (or a new set of medians) could be computed.  Even if the computational effort is not large, an excessive number of solutions may be unacceptable for human users. Instead, with this approach a set of solutions is computed once and then we can choose the best one in real time taken from a relatively small set of solutions.\\

The parameters of the uncertainty set  may change over the time and the scenarios appear according it. In that case a solution with robustness close to $v(R(\Gamma))$ will be available any time for the {\it true} $\Gamma$. In practice, neither the true parameters nor a corresponding robust solution need to be known for the decision maker at any time. Obviously we need to estimate $\Omega$.\\

The approach to be presented is complementary to the min-max-min approach to find $k$ solutions that work well for a fixed $\Gamma$ (\cite{BuchKurtz2016},\cite{BuchKurtz2017},\cite{Kurtz2021}).\\

Formally: let ${\lbrace x^i \rbrace}_1^r \subseteq X$ and let $\epsilon \geq 0$. We say that {\it ${\lbrace x^i \rbrace}_1^r$ is an $\epsilon$-optimal multiparametric robust solution for $\Omega$} ($\epsilon,\Omega$-mprs) if  
 
\[\underset{i \in [r]}{\min}\;v(W( x^i ,\Gamma))-v(R(\Gamma)) \leq \epsilon\; \forall \Gamma \in \Omega\]

Some problems will be presented several times with equivalent formulations. In some cases the name of the problem will be the same and the specific formulation used would be clear from the context. In order to clarify the exposition in some cases we will use different names for the same problem according the formulation used.\\ 

The paper is organized as follows: in section 2 we present an algorithm to  find an $\epsilon,\Omega$-mprs. The algorithm consists of applying an algorithm to obtain an $\epsilon,\Omega$-{\it optimal} multiparametric solution relative to the objective function (\cite{Cre2000},\cite{Cre2020}) for a 0-1-Integer Linear Programming (0-1-ILP) problem  equivalent to $R(\Gamma)$. The case in which the matrix that represents $X$ is totally unimodular is considered in a subsection. In section 3 we present the details for the cases by interval, linear segment and budgeted. Computational experience is presented in section 4 for Shortest Path and $(l,p)$-Medians problems by using ILOG-Cplex 12.10 from a DOcplex Python code (\cite{Cplex}) . In section 5  we present a 0-1-MILP problem equivalent to $R(\Gamma)$ for the variant  for the uncertainty set defined as follows (\cite{GoeLen2021}):
 
\[\Lambda^+(\Gamma) = \lbrace c \in \mathbb{R}^n_+: c_j = \underline{c}_j + \lambda_j d_j\;\forall j \in [n],\;\underset{j \in P_k}{\sum} \lambda_j \leq \Gamma_k\;\forall k \in [K],\;\lambda \in [0,1]^n\rbrace\]

in such a manner that a multiparametric analysis may be performed to find an $\epsilon,\Omega$-mprs.    Conclusions and further extensions may be seen in section 6.\\

Appendices may be seen after references. In Appendix A we present a summary of a multiparametric approach relative to the objective function for 0-1-ILP problems. In Appendix B we present a toy example to show that with low uncertainty a large number of solutions may be necessary to define a $0,\Omega$-mprs and with large uncertainty a single solution may be enough to define a $0,\Omega$-mprs. In Appendix C we present a proof of properties for an auxiliary function that appears in section 5. A remark about the computational complexity is presented in Appendix D. Tables may be seen after appendices.

\section{Algorithm to find an $\epsilon,\Omega$-multiparametric robust solution}

Let $x \in X$ and let $\Gamma \in \Omega$, the dual problem of $W(x,\Gamma)$ is defined as a LP problem in $(\pi,\rho)$ as follows
\begin{flalign*}
\underline{c}^t x + \min&\;\Gamma^t \pi +  d^t \rho&DW(x,\Gamma)&\\
s.t.&\;\;\pi_k + \rho_j \geq x_j&\forall k \in [k]\;\forall j \in P_k&\\
&\;\;\pi \geq 0,\;\rho \geq 0&\\
&\;\;\pi \in \mathbb{R}^K,\;\rho \in \mathbb{R}^n&\
\end{flalign*}

There exists an optimal solution for $DW(x,\Gamma)$ with $\pi \in {\lbrace 0,1 \rbrace}^K$ (\cite{GoeLen2021}), therefore there exists an optimal solution for $DW(x,\Gamma)$ with $(\pi,\rho) \in {\lbrace 0,1 \rbrace}^K \times {\lbrace 0,1 \rbrace}^n$, hence $R(\Gamma)$ may be rewritten as a 0-1-ILP problem in $(\pi,\rho,x)$ as follows:
\begin{flalign*}
\min&\;\;\Gamma^t \pi + d^t \rho + \underline{c}^t  x&R(\Gamma)&\\
s.t.&\;\;\pi_k + \rho_j - x_j \geq 0&\forall k \in [k]\;\forall j \in P_k&\\
&\;\;\pi \in {\lbrace 0,1 \rbrace}^K,\;\rho \in {\lbrace 0,1 \rbrace}^n,\;x \in X&\
\end{flalign*}

Let $\epsilon \geq 0$.  We may use a multiparametric algorithm (see Appendix A) to find an $\epsilon,\Omega$-optimal multiparametric solution for $R(\Gamma)$, that is we may find 
$\lbrace ({\pi^i,\rho^i,x^i) \rbrace}_1^r \subseteq {\lbrace 0,1 \rbrace}^K \times {\lbrace 0,1 \rbrace}^n \times X$ such that: $\pi^i_k + \rho^i_j - x^i_j \geq 0\;\forall k \in [k]\;\forall j \in P_k,\;\forall i \in [r]$ and $\underset{i \in [r]}{\min}\;\left \lbrace \Gamma^t \pi^i + d^t \rho^i + \underline{c}^t  x^i \right \rbrace - v(R(\Gamma)) \leq \epsilon\;\;\forall \Gamma \in \Omega$. In that case because of $v(W(x^i,\Gamma)) = v(DW(x^i,\Gamma)) \leq \Gamma^t \pi^i + d^t \rho^i + \underline{c}^t  x^i$ for all $i \in [r]$ we have that ${\lbrace x^i \rbrace}_1^r$ is an $\epsilon,\Omega$-mprs\\

Next we present the multiparametric algorithm following Appendix A applied to $R(\Gamma)$.\\ 

Let ${\lbrace (\pi^i,\rho^i,x^i) \rbrace}_1^r \subseteq {\lbrace 0,1 \rbrace}^K \times {\lbrace 0,1 \rbrace}^n \times X$ and let $Q({\lbrace (\pi^i,\rho^i,x^i) \rbrace}_1^r)$ be a problem in $(\Gamma,\pi,\rho,x)$ defined as follows: 
\begin{flalign*}
\max& \left (\underset{i \in [r]}{\min} \left \lbrace \Gamma^t  \pi^i + d^t \rho^i + \underline{c}^t x^i \right \rbrace- \left ( \Gamma^t \pi +  d^t \rho + \underline{c}^t x \right ) \right )&Q({\lbrace (\pi^i,\rho^i,x^i) \rbrace}_1^r)&\\
s.t.&\;\;\pi_k + \rho_j - x_j \geq 0&\forall k \in [k]\;\forall j \in P_k&\\
&\;\;\Gamma \in \Omega,\;\pi \in {\lbrace 0,1 \rbrace}^K,\;\rho \in {\lbrace 0,1 \rbrace}^n,\;x \in X&\
\end{flalign*}

{\bf Algorithm $Q$ (A-$Q$) to find an $\epsilon,\Omega$-multiparametric robust solution}\\

Let $\epsilon \geq 0$ and let $\Gamma \in \Omega$. Solve $R(\Gamma)$, let $(\pi^1,\rho^1,x^1)$ be an optimal solution and let $r=1$. 

\begin{enumerate}  
	\item Solve $Q({\lbrace (\pi^i,\rho^i,x^i) \rbrace}_1^r)$ and let $(\Gamma^*,\pi^*,\rho^*,x^*)$ be an optimal solution
	\item If $v(Q({\lbrace (\pi^i,\rho^i,x^i) \rbrace}_1^r)) \leq \epsilon$ STOP
	\item Let $(\pi^{r+1},\rho^{r+1},x^{r+1}) =(\pi^*,\rho^*,x^*) $, let $r = r+1$ and return to step 1
\end{enumerate}

If $(\Gamma,\pi,\rho,x)$ is an optimal solution for $Q({\lbrace (\pi^i,\rho^i,x^i) \rbrace}_1^r)$ then $(\pi,\rho,x)$ is an optimal solution for $R(\Gamma)$. Therefore with $Q({\lbrace (\pi^i,\rho^i,x^i) \rbrace}_1^r)$ we are looking for $\Gamma \in \Omega$ that maximizes $\underset{i \in [r]}{\min} \left \lbrace \Gamma^t  \pi^i + d^t \rho^i + \underline{c}^t x^i \right \rbrace - v(R(\Gamma))$. If the difference is less or equal to $\epsilon$ we are done. Otherwise we updated ${\lbrace (\pi^i,\rho^i,x^i) \rbrace}_1^r$ by introducing $(\pi,\rho,x)$. Since ${\lbrace 0,1 \rbrace}^K \times {\lbrace 0,1 \rbrace}^n \times X$ is a finite set {\bf A-$Q$} stops in a finite number of iterations. If ${\lbrace (\pi^i,\rho^i,x^i) \rbrace}_1^r$ is the output then ${\lbrace x^i \rbrace}_1^r$ is an $\epsilon,\Omega$-mprs.\\ 

In order to solve $Q({\lbrace (\pi^i,\rho^i,x^i) \rbrace}_1^r)$ we can rewrite it as a 0-1-MILP problem in $(\Gamma,\pi,\rho,x,w,\sigma)$ as follows:
\begin{flalign*}
\max&\;\;\sigma - \left ( \underset{k \in [K]}{\sum} w_k +  d^t \rho + \underline{c}^t x \right )&Q({\lbrace (\pi^i,\rho^i,x^i) \rbrace}_1^r)&\\
s.t.&\;\;\sigma - \Gamma^t \pi^i \leq  d^t \rho^i + \underline{c}^t x^i&\forall i \in [r]&\\
&\;\;\pi_k + \rho_j - x_j \geq 0&\forall k \in [k]\;\forall j \in P_k&\\
&\;\;w_k - \Gamma_k - \mathcal{U}_k \pi_k \geq - \mathcal{U}_k&\forall k \in [K]&\\
&\;\;\Gamma \in \Omega,\;\pi \in {\lbrace 0,1 \rbrace}^K,\;\rho \in {\lbrace 0,1 \rbrace}^n,\;x \in X&\\
&\;\;w \in \mathbb{R}^K_+,\;\sigma \in \mathbb{R}&\
\end{flalign*}

Let $(\Gamma,\pi,\rho,x,w,\sigma)$ be an optimal solution for $Q({\lbrace (\pi^i,\rho^i,x^i) \rbrace}_1^r)$. If $\pi_k = 1$ then $w_k \geq \Gamma_k$ and because of the maximization criterium we have $w_k = \Gamma_k = \Gamma_k \pi_k$. If $\pi_k = 0$ then $w_k \geq \Gamma_k - \mathcal{U}_k$ and because of the maximization criterium
we have $w_k = 0 = \Gamma_k \pi_k$. Therefore we have $w_k = \Gamma_k \pi_k$.  Since maximization is the criterium we have that $\sigma = \underset{i \in [r]}{\min}\;\left \lbrace  \Gamma^t \pi^i + d^t \rho^i + \underline{c}^t x^i \right \rbrace$. Therefore, the reformulation of the $Q$ problem is valid.  

\subsection{Totally unimodular case}

Let $c(\pi)_j = (\underline{c}_j \pi_k) +  (\underline{c}_j+d_j)(1-\pi_k)\;\forall j\in P_k\;\forall k \in [K]$. After some algebraic manipulations
we have that $R(\Gamma)$ may be rewritten as follows (\cite{GoeLen2021}):
\begin{flalign*}
\min&\;\;\Gamma^t \pi + v(P(c(\pi)))&\hat R(\Gamma)&\\
s.t.&\;\;\pi \in {\lbrace 0,1 \rbrace}^K&\
\end{flalign*}

and we know that if $(\pi,\rho,x)$ is an optimal solution for $R(\Gamma)$ then $\pi$ is an optimal solution for $\hat R(\Gamma)$ and $x$ is an optimal solution for $P(c(\pi))$. By the other hand if (i) $\pi$ is an optimal solution for $\hat R(\Gamma)$, (ii) $x$ is an optimal solution for $P(c(\pi))$ and (iii) $\rho_j = (1-\pi_k) x_j\;\forall j \in P_k\;\forall k \in [K]$ then $(\pi,\rho,x)$ is an optimal solution for $R(\Gamma)$.\\

Let $\overline{P}(c)$ be the linear relaxation of $P(c)$ with $\overline{X}$ instead of $X$ by using $x \in [0,1]^n$ instead of $x \in {\lbrace 0,1 \rbrace}^n$. Let $\overline{R}(\Gamma)$ be a 0-1-MILP problem in $(\pi,\rho,x)$ defined from $R(\Gamma)$ with  $(\rho,x) \in [0,1]^n \times [0,1]^n$ instead of $(\rho,x) \in {\lbrace 0,1 \rbrace}^n \times {\lbrace 0,1 \rbrace} ^n$. We have that  $\overline{R}(\Gamma)$  may be rewritten as follows:
\begin{flalign*}
\min&\;\;\Gamma^t \pi + v(\overline{P}(c(\pi)))&\hat{\overline{R}}(\Gamma)&\\
s.t.&\;\;\pi \in {\lbrace 0,1 \rbrace}^K&\
\end{flalign*}

and we know that if $(\pi^*,\bar \rho,\bar x)$ is an optimal solution for $\overline{R}(\Gamma)$ then $\pi^*$ is an optimal solution for $\hat{ \overline{R}}(\Gamma)$ and $\bar x$ is an optimal solution for $\overline{P}(c(\pi^*))$. By the other hand if (i) $\pi^*$ is an optimal solution for $\hat {\overline{R}}(\Gamma)$, (ii) $\bar x$ is an optimal solution for $\overline{P}(c(\pi^*))$ and (iii) $\bar \rho_j = (1-\pi^*_k) \bar x_j\;\forall j \in P_k\;\forall k \in [K]$ then $(\pi^*,\bar \rho, \bar x)$ is an optimal solution for $\overline{R}(\Gamma)$.\\

Let us suppose that $P(c)$ is defined as follows:
\begin{flalign*}
\min&\;\; c^t x&P(c)&\\
s.t.&\;\;Ax = b&\\
&\;\;x \in {\lbrace 0,1 \rbrace}^n&\
\end{flalign*}

with $b \in \mathbb{Z}^m$ and $A \in \mathbb{R}^{m \times n}$.\\

Let us suppose that $A$ is totally unimodular. We know that in the totally unimodular case we have that $v(P(c)) = v(\overline{P}(c))$ for all $c$ and if $\bar x$ is an optimal solution for $\overline{P}(c)$ then $\bar x \in {\lbrace 0,1 \rbrace}^n$ and $\bar x$ is an optimal solution for $P(c)$. Therefore, in the totally unimodular case if $(\pi^*,\bar \rho,\bar x)$ is an optimal solution for $\overline{R}(\Gamma)$ we have that $\bar x$ is an optimal solution for $\overline{P}(c)$ with $\bar x \in {\lbrace 0,1 \rbrace}^n$. If $\rho^*_j = (1-\pi^*_k) \bar x_j$ for all $j \in P_k$ and for all $k \in [K]$ then $\rho^* \in {\lbrace 0,1 \rbrace}^n$, therefore :  $(\pi^*,\rho^*, \bar x)$ is an optimal solution for $R(\Gamma)$\\ 

Let $\overline{Q}({\lbrace (\pi^i,\rho^i,x^i) \rbrace}_1^r)$ be a 0-1-MILP problem in $(\Gamma,\pi,\rho,x,w,\sigma)$ defined from $Q({\lbrace (\pi^i,\rho^i,x^i) \rbrace}_1^r)$ with $(\rho,x) \in [0,1]^n \times \overline{X}$ instead of $(\rho,x) \in {\lbrace 0,1 \rbrace}^n \times X$.\\

If $(\Gamma,\pi,\rho,x,w,\sigma)$ is an optimal solution for $\overline{Q}({\lbrace (\pi^i,\rho^i,x^i) \rbrace}_1^r)$ then \break  $(\pi,\rho,x)$ is an optimal solution for $\overline{R}(\Gamma)$. Since we have the totally unimodular case $(\pi,\rho^*,x)$ is an optimal solution for $R(\Gamma)$ with $\rho^*_j =(1-\pi_k)  x_j$ for all $j \in P_k$ and for all $k \in [K]$, hence $(\Gamma,\pi,\rho^*,x,w,\sigma)$ is an optimal solution for $Q({\lbrace (\pi^i,\rho^i,x^i) \rbrace}_1^r)$.

\section{Specific cases for $\Omega$}

Next we present the $Q$ problem formulation for three $\Omega$-cases of practical interest. We use standard formulation procedures and some algebraic manipulations. 

\subsection{Interval case}

Let $\mathcal{L} \in \mathbb{R}^K_+$ with $\mathcal{L} \leq \mathcal{U}$ and let $\Omega = [\mathcal{L},\mathcal{U}] = \lbrace \Gamma \in \mathbb{R}^K_+: \mathcal{L}  \leq \Gamma \leq \mathcal{U} \rbrace$\\

Let $\Gamma^+(\pi)_k = \mathcal{L}_k \pi_k + \mathcal{U}_k (1-\pi_k)$ for all $k \in [K]$, for all $\pi \in {\lbrace 0,1 \rbrace}^K$.\\

According to appendix A:  $Q({\lbrace (\pi^i,\rho^i,x^i) \rbrace}_1^r)$ 
is equivalent to $Q^+({\lbrace (\pi^i,\rho^i,x^i) \rbrace}_1^r)$ defined as follows:
\begin{flalign*}
\max&\;\;\underset{i \in [r]}{\min} \left \lbrace {\Gamma^+(\pi)}^t  \pi^i + d^t \rho^i + \underline{c}^t x^i \right \rbrace- \left ( \mathcal{L}^t \pi +  d^t \rho + \underline{c}^t x \right ) &Q^+({\lbrace (\pi^i,\rho^i,x^i) \rbrace}_1^r)&\\
s.t.&\;\;\pi_k + \rho_j - x_j \geq 0&\forall k \in [k]\;\forall j \in P_k&\\
&\;\;\pi \in {\lbrace 0,1 \rbrace}^K,\;\rho \in {\lbrace 0,1 \rbrace}^n,\;x \in X&\
\end{flalign*}

In order to solve $Q^+({\lbrace (\pi^i,\rho^i,x^i) \rbrace}_1^r)$ we may rewrite it as a 0-1-ILP in $(\pi,\rho,x,\sigma)$ as follows:
\begin{flalign*}
\max&\;\;\sigma - \left (\mathcal{L}^t \pi +  d^t \rho + \underline{c}^t x \right )&Q^+({\lbrace (\pi^i,\rho^i,x^i) \rbrace}_1^r)&\\
s.t.&\;\;\sigma + {f^i}^t \pi \leq  d^t \rho^i + \underline{c}^t x^i + \mathcal{U}^t \pi^i&\forall i \in [r]&\\
&\;\;\pi_k + \rho_j - x_j \geq 0&\forall k \in [k]\;\forall j \in P_k&\\
&\;\;\pi \in {\lbrace 0,1 \rbrace}^K,\;\rho \in {\lbrace 0,1 \rbrace}^n,\;x \in X,\;\sigma \in \mathbb{R}&\
\end{flalign*}

with $f^i_k = (\mathcal{U}_k -\mathcal{L}_k) \pi^i_k\;\forall k \in [K],\;\forall i \in [r]$\\

Since maximization is the criterium we have that\\
$\sigma = \underset{i \in [r]}{\min}\;\left \lbrace d^t \rho^i + \underline{c}^t x^i + \mathcal{U}^t \pi^i - {f^i}^t \pi\right \rbrace  = \underset{i \in [r]}{\min}\;\left \lbrace {\Gamma^+(\pi)}^t \pi^i + d^t \rho^i + \underline{c}^t x^i \right \rbrace$ for all optimal solution. Therefore, the reformulation of the $Q^+$ problem is valid.

\subsection{Line segment case}

Let $\Gamma^0 \in \mathbb{R}^K_+$, let $[\underline{\alpha},\overline{\alpha}] \subseteq [0,1]$ and let $\Omega = \lbrace \Gamma \in \mathbb{R}^K_+:\Gamma = \alpha  \Gamma^0,\;\alpha \in [\underline{\alpha},\overline{\alpha}] \rbrace$ with $\mathcal{U} = \overline{\alpha}  \Gamma^0$\\

$Q({\lbrace (\pi^i,\rho^i,x^i) \rbrace}_1^r)$ becomes a 0-1-MILP problem in $(\alpha,\pi,\rho,x,w,\sigma)$ defined as follows:
\begin{flalign*}
\max&\;\;\sigma - \left ( \underset{k \in [K]}{\sum} w_k +  d^t \rho + \underline{c}^t x \right )&Q({\lbrace (\pi^i,\rho^i,x^i) \rbrace}_1^r)&\\
s.t.&\;\;\sigma - \alpha {\Gamma^0}^t \pi^i \leq  d^t \rho^i + \underline{c}^t x^i&\forall i \in [r]&\\
&\;\;\pi_k + \rho_j - x_j \geq 0&\forall k \in [K]\;\forall j \in P_k&\\
&\;\;w_k - \alpha \Gamma^0_k - \overline{\alpha} \Gamma^0_k \pi_k \geq -\overline{\alpha}\Gamma^0_k &\forall k \in [K]&\\
&\;\;\alpha \in [\underline{\alpha},\overline{\alpha}],\;\pi \in {\lbrace 0,1 \rbrace}^K,\;\rho \in {\lbrace 0,1 \rbrace}^n,\;x \in X&\\
&\;\;w \in \mathbb{R}^K_+,\;\sigma \in \mathbb{R}&\
\end{flalign*}

Let $(\alpha,\pi,\rho,x,w,\sigma)$ be an optimal solution for $Q({\lbrace (\pi^i,\rho^i,x^i) \rbrace}_1^r)$. If $\pi_k = 1$ then $w_k - \alpha \Gamma^0_k \geq 0$ and because of the maximization criterium we have $w_k = \alpha \Gamma^0_k = \alpha \Gamma^0_k \pi_k = \Gamma_k \pi_k$. If $\pi_k = 0$ then $w_k \geq (\alpha - \overline{\alpha}) \Gamma^0_k$ and because of the maximization criterium
we have $w_k = 0 = \Gamma_k \pi_k$. Therefore we have $w_k = \Gamma_k \pi_k$ for all $k \in [K]$. Since maximization is the criterium we have that\\

$\sigma = \underset{i \in [r]}{\min}\;\left \lbrace  \alpha {\Gamma^0}^t \pi^i + d^t \rho^i + \underline{c}^t x^i \right \rbrace  = \underset{i \in [r]}{\min}\;\left \lbrace \Gamma^t \pi^i + d^t \rho^i + \underline{c}^t x^i \right \rbrace$.\\
		
Therefore, the reformulation of the $Q$ problem is valid.

\subsection{Budgeted case}

Let $\underline{\Gamma},\mathcal{D} \in \mathbb{R}^K_+$, let $\Delta \in \mathbb{R}_+$ and let $\Omega =\lbrace \Gamma \in \mathbb{R}^K_+:\Gamma = \underline{\Gamma} + \beta,\;\;\underset{k \in [K]}{\sum} \beta_k \leq \Delta,\;\beta_k \in [0,\mathcal{D}_k]\;\forall k \in [K] \rbrace$ with  $\mathcal{U} = \underline{\Gamma} + \mathcal{D}$\\

$Q({\lbrace (\pi^i,\rho^i,x^i) \rbrace}_1^r)$ becomes a 0-1-MILP problem in $(\beta,\pi,\rho,x,w,\sigma)$ defined as follows:
\begin{flalign*}
\max&\;\;\sigma - \left ( \underset{k \in [K]}{\sum} w_k +  d^t \rho + \underline{c}^t x \right )&Q({\lbrace (\pi^i,\rho^i,x^i) \rbrace}_1^r)&\\
s.t.&\;\;\sigma - \beta^t \pi^i \leq  d^t \rho^i + \underline{c}^t x^i + \underline{\Gamma}^t \pi^i &\forall i \in [r]&\\
&\;\;\pi_k + \rho_j - x_j \geq 0&\forall k \in [K]\;\forall j \in P_k&\\
&\;\;w_k - \beta_k -  (\underline{\Gamma} + \mathcal{D})_k \pi_k \geq -  \mathcal{D}_k &\forall k \in [K]&\\
&\;\;\beta_k \leq \mathcal{D}_k&\forall k \in [K]&\\
&\;\;\underset{k \in [K]}{\sum} \beta_k \leq \Delta&\\
&\;\;\beta \in \mathbb{R}^K_+,\;\pi \in {\lbrace 0,1 \rbrace}^K,\;\rho \in {\lbrace 0,1 \rbrace}^n,\;x \in X&\\
&\;\;w \in \mathbb{R}^K_+,\;\sigma \in \mathbb{R}&\
\end{flalign*}

Let $(\beta,\pi,\rho,x,w,\sigma)$ be an optimal solution for $Q({\lbrace (\pi^i,\rho^i,x^i) \rbrace}_1^r)$. If $\pi_k = 1$ then $w_k \geq \underline{\Gamma}_k + \beta_k$ and because of the maximization criterium we have $w_k = \underline{\Gamma}_k + \beta_k = (\underline{\Gamma}_k + \beta_k) \pi_k$ for all $k \in [K]$. If $\pi_k = 0$ then $w_k \geq \beta_k - \mathcal{D}_k$  and because of the maximization criterium
we have $w_k = 0 =  (\underline{\Gamma}_k + \beta_k)\pi_k$. Therefore we have $w_k = \Gamma_k \pi_k$ for all $k \in [K]$.  Since maximization is the criterium we have that\\

$\sigma = \underset{i \in [r]}{\min}\;\left \lbrace (\underline{\Gamma} + \beta)^t \pi^i  + d^t \rho^i + \underline{c}^t x^i \right \rbrace  = \underset{i \in [r]}{\min}\;\left \lbrace \Gamma^t \pi^i + d^t \rho^i + \underline{c}^t x^i \right \rbrace$.\\

Therefore, the reformulation of the $Q$ problem is valid.

\section{Computational experience} 

Our algorithms have been performed on a personal computer as follows:
\begin{itemize}
	\item Intel(R)Core(TM) i7-9750H CPU, @ 2.60 GHz Lenovo ThinkPad X1 Extreme Gen 2, 32.00 GB Ram and Windows 10 Pro Operating System
	\item All the instances have been processed through ILOG-Cplex 12.10 
	from a DOcplex Python code
	\item  All the parameters of ILOG-Cplex 12.10 are in their default values
\end{itemize}

\subsection{Shortest path problem}

\subsubsection{Data generation}

Graphs $G = (V,E)$ are  generated following \cite{Han2015} as follows: in each problem instance, the nodes correspond to $|V|$ points with coordinates  that are chosen uniformly at random in the square $[0,10] \times [0,10]$. We choose the pair of nodes with the largest Euclidean distance as $v_0$ and $w_0$ (the start and terminal nodes). In order to generate $E$, we begin with a fully connected graph and remove $70\%$ of the arcs  in order of decreasing Euclidean distance, that is, starting with the longest arcs. If $e=(v,w) \in E$ then $\underline{c}_e$ is the euclidean distance from $v$ to $w$ and $d_e = 0.5 \underline{c}_e$.\\
                                 
Let $b_{v_0} = 1,\;b_{w_0} = 1$ and $b_v  = 0\;\forall v \in V - \lbrace v_0,w_0 \rbrace$.
For each $v \in V$ let $E^+(v) = \lbrace (v,w): (v,w) \in E  \rbrace$ and let $E^-(v) = \lbrace (w,v):(w,v) \in E \rbrace$.\\
                                 
The shortest path (SP) problem for $c$ is defined as follows: 
\begin{flalign*}                 
\min&\;\;\underset{e \in E}{\sum} c_e x_e&SP(c)&\\
s.t.&\sum_{e \in E^+(v)} x_e  - \sum_{e \in E^-(v)} x_e =b_v&\forall v \in V&\\
&\;\;x \in {\lbrace 0,1 \rbrace}^{|E|}&\
\end{flalign*}

\subsubsection{Partitions definitions}

\begin{itemize}
\item Random generation of the partitions (marked with $r$): for each $e \in E$ an index is randomly selected in $[K]$, if the selected index is $k$ then $e \in P_k$  

\item Generation of the partitions according values of a minimal path to destination (marked with $p$): let $cost(v)$ the cost of a minimum path from $v$ to $w_0$ for all $v \in V$. Let $\hat {cost} = \underset{v \in V}{\max}\; cost(v)$ and let $\lambda = \hat {cost}/K$. If $cost(v) \in  [(k-1) \lambda,k \lambda]$ then $e \in P_k$ for all $e \in \delta^+(v)$   

\item Generation of the partitions according distances to destination (marked with $d$): let $dist(v)$ the euclidean distance from $v$ to $w_0$ for all $v \in V$. Let $\hat {dist} = \underset{v \in V}{\max}\; dist(v)$ and let $\lambda = \hat {dist}/K$. If $dist(v) \in  [(k-1)\lambda,k \lambda]$ then $e \in P_k$ for all $e \in \delta^+(v)$
\end{itemize}

\subsubsection{$\Omega$ definition:}

\begin{itemize}
	\item Interval case: let $[\mathcal{L}_k(\delta),\mathcal{U}_k(\delta)] = \underset{e \in P_k}{\max} \lbrace d_e \rbrace \times [\delta,\delta+1]\;\;\forall k \in [K]$ with $\delta \geq 0$ a parameter to be presented in the tables. Note that the volumes of the $\Omega$ sets used are equals ($\underset{k \in [K]}\Pi (\mathcal{U}_k(\delta) -\mathcal{L}_k(\delta)) = \underset{k \in [K]}\Pi \underset{e \in P_k}{\max} \lbrace d_e \rbrace\;\forall \delta \geq 0$) 
	\item Line segment case: $\Gamma^0_k = \underset{e \in P_k}{\max} \lbrace d_e \rbrace\;\;\forall k \in [K],\;\alpha \in  [\underline{\alpha},\overline{\alpha}] = [0,1]$
	\item Budgeted case: $\underline{\Gamma}_k = \beta_1 \times  \underset{e \in P_k}{\max} \lbrace d_e \rbrace\;\;\forall k \in [K]$, $\mathcal{D}=\beta_2 \underline{\Gamma}$, $\Delta = \delta \underset{k \in [K]}{\max} \lbrace \mathcal{D}_k \rbrace$ with $\beta_1 \geq 0,\beta_2 \geq 0,\delta \geq 0$   parameters to be presented in the tables    
\end{itemize}

\subsection{(l,p)-Medians problem}

\subsubsection{Data generation}

Let $[l]$ a set of demand locations.  Each demand location is a candidate to be a service location (a median). If the demand location $j$ is assigned to the median $i$ the cost is $c_{ij}$ . Let $p$ the number of medians to be selected.\\

The data were generated at random as follows: locations $j$ were taken from $U((0,100) \times (0,100))$, let $D_j$ be the demand of location $j$ taken from $U(0,100)$, let $dist_{ij}$ be the  euclidean distance from location $i$ until location $j$ and let $\underline{c}_{ij} = dist_{ij} D_j\;\;\forall i,j \in [l]$. If $(i,j) \in [l] \times [l]$ then $d_{ij} = 0.5 \underline{c}_{ij}$.\\

Let $y_i \in \lbrace 0,1 \rbrace\;\;(i \in [l])$ with $y_i = 1$ if and only if the demand location $i$ is selected to be a median. Let $x_{ij} \in \lbrace 0,1 \rbrace \;\;(i,j) \in [l] \times [l])$
with $x_{ij} = 1$ if and only if the demand location $j$ is assigned to a median located at $i$.\\ 

The $(l,p)$-medians ($(l,p)$M) problem with cost $c$ is a problem in $(y,x)$ defined as follows (\cite{BofKar1984}):
\begin{flalign*}
\min& \sum_{i\in [l]} \sum_{j \in [l]} c_{ij} x_{ij}&(l,p)M(c)\\
s.t.&\;\;x_{ij} \leq y_i&\forall i,j \in [l]&\\
&\;\;\sum_{i \in [l]} y_i = p&\\
&\;\;\sum_{i \in [l]} x_{ij} = 1&\forall j \in [l]&\\
&\;\;y_i \in \lbrace 0,1 \rbrace,\;\;x_{ij} \in \lbrace 0,1 \rbrace&\forall i,j \in [l]&\
\end{flalign*}
	
\subsubsection{Partitions definition}

\begin{itemize}
	\item Generation of the partitions guided by locations (marked with $lo$): for each $i \in [l]$ an index is randomly selected in $[K]$, if the selected index is $k$ then $(i,j) \in P_k$ for all $j \in [l]$
	\item Generation of the partitions according the sum of perturbations (marked with $g$): let $cost(i) = \underset{j \in [l]}{\sum} d_{ij}$. Let $\hat {cost} = \underset{i \in [l]}{\max}\; cost(i)$ and let $\lambda = \hat {cost}/K$. If $cost(i) \in  [(k-1) \lambda,k \lambda]$ then $(i,j) \in P_k$ for all $j \in [l]$         
\end{itemize}

\subsubsection{$\Omega$ definition:}

\begin{itemize}
	\item Interval case: $[\mathcal{L}_k(\delta),\mathcal{U}_k(\delta)] = \underset{(i,j) \in P_k}{\max} \lbrace d_{ij} \rbrace \times [\delta,\delta+1]\;\;\forall k \in [K]$ with $\delta \geq 0$ a parameter to be presented in the tables. Note that the volumes of the $\Omega$ sets used are equals ($\underset{k \in [K]}\Pi (\mathcal{U}_k(\delta) -\mathcal{L}_k(\delta)) = \underset{k \in [K]}\Pi \underset{(i,j) \in P_k}{\max} \lbrace d_{ij} \rbrace\;\forall \delta \geq 0$)
	\item Line segment case: $\Gamma^0_k = n \times \underset{(i,j) \in P_k}{\max} \lbrace d_{ij} \rbrace\;\;\forall k \in [K],\;\alpha \in  [\underline{\alpha},\overline{\alpha}] = [0,1]$
	\item Budgeted case: $\underline{\Gamma}_k = \beta_1 \times  \underset{(i,j) \in P_k}{\max} \lbrace d_{ij} \rbrace\;\;\forall k \in [K]$, $\mathcal{D}=\beta_2 \underline{\Gamma}$, $\Delta = \delta \underset{k \in [K]}{\max} \lbrace \mathcal{D}_k \rbrace$ with $\delta \geq 0,\beta_1 \geq 0,\beta_2 \geq 0$  parameters to be presented in the tables.    
\end{itemize}

\subsection{Conditions to start and stop the algorithm}

\begin{itemize}
	\item Interval case: solve $R(\mathcal{L})$ and let $(x^1,\pi^1,\rho^1)$ be an optimal solution 
	\item Budgeted case: solve $R(\underline{\Gamma})$ and let $(x^1,\pi^1,\rho^1)$ be an optimal solution 
	\item Line segment case: solve $R(0)$ and let $(x^1,\pi^1,\rho^1)$ be an optimal solution
\end{itemize}

We use $\epsilon$ as follows:

\begin{itemize}
	\item Interval case: $\epsilon = 0.01 v(R(\mathcal{L}))$ 
	\item Budgeted case: $\epsilon = 0.01 v(R(\underline{\Gamma}))$
	\item Line segment case: $\epsilon = 0.01 v(R(0)) = 0.01 v(P(\underline{c}))$
\end{itemize}

Note that if ${\lbrace x^i \rbrace}_1^r$ is the output generated by the algorithm then we have:
\[v(R(\Gamma)) \leq \underset{i \in [r]}{\min}\;v(W( x^i ,\Gamma)) \leq (1+0.01) v( R(\Gamma))\; \forall \Gamma \in \Omega\]
Note that in each iteration $v(Q({\lbrace \pi^i,\rho^i,x^i \rbrace}_1^r))/v(R(\mathcal{L}))$ is an upper bound of the maximal relative error defined as $\underset{\Gamma \in \Omega}{\max}\;\frac{\underset{i \in [r]}{\min}\;v(W( x^i ,\Gamma)) -v(R(\Gamma))}{v(R(\Gamma))}$ 

\subsection{Values to be presented}

Each line in the tables refer to a set of problems under conditions presented  and we present:

\begin{itemize}
	\item $\overline{t_*}$ and $\hat {t_*}$: the average and maximun time in seconds (s.) to find an $\epsilon,\Omega$-mprs when the partitions are generated according $* \in \lbrace r,p,d,lo,g \rbrace$, including the time to find a robust solution ($x^1$) 
	\item $\underline{s_*},\;\overline{s_*}$ and $\hat s_*$: the minimum, average and maximun number of different solutions (paths to SP problem and $p$-medians for $(l,p)$M-problem) to define an $\epsilon,\Omega$-mprs  when  the partitions are generated according  $* \in \lbrace r,p,d,lo,g \rbrace$. Do not confuse these values with the corresponding minimum, average and maximum number of iterations executed by the algorithm, which can be considerably higher in some cases. The reason is that if $(\pi^1,\rho^1,x^1)$ and $(\pi^2,\rho^2,x^2)$ are different solutions generated by the algorithm, it often happens that $x^1 = x^2$. The decision maker is only interested in bailing out $x=x^1=x^2$ in cases like that.\\ 
	
	For $(l,p)$M-problems we may find a lot of cases with different solutions with the same $p$-medians set. For the decision maker the hard decision is to choice a $p$-medians set because of the assignment of locations to medians is very simple in real time. Hence, we report minimum, average and maximun number of different $p$-medians set to define an $\epsilon,\Omega$-mprs   
	\item $\underline{s_{d,q}}$, $\overline{s_{d,q}}$ and $\hat{s_{d,q}}$: the minimun, average and maximun number of different solutions (paths) to define a $\epsilon(q),\Omega$-mprs with $\epsilon(q) = (q/100) \times  v(R(\mathcal{L}))$ when the $P_k$ definition corresponds to $d$ for $q \in \lbrace 1,2,3,5 \rbrace$
	\item  $\overline{t_{d,q}}$: the average time in seconds to find $\epsilon(q),\Omega$-mprs with $\epsilon(q) = (q/100) \times v(R(\mathcal{L}))$ when the $P_k$ definition corresponds  to  $d$ for $q \in \lbrace 1,2,3,5 \rbrace$, including the time to find a robust solution ($x^1$)
	\item We present, according the case:
	\begin{enumerate}
	  \item $\overline{\epsilon_{\mathcal{L}}}$: the average of $ 100 \times \frac{v(Q(\pi^1,\rho^1,x^1))}{v(R(\mathcal{L}))}$ 
	  \item $\overline{\epsilon_{\underline{\Gamma}}}$: the average of $100 \times \frac{v(Q(\pi^1,\rho^1,x^1))}{v(R(\underline{\Gamma}))}$ 
	  \item $\overline{ \epsilon_0}$: the average of $100 \times  \frac{v(Q(\pi^1,\rho^1,x^1))}{v(R(0))}$ 
	  \item $\overline{ \epsilon^{(2)}_0}$: the average of $100 \times  \frac{v(Q({\lbrace \pi^i,\rho^i,x^i \rbrace}_1^2 ))}{v(R(0))}$ 
	  
	  \end{enumerate}
\end{itemize}

All values were rounded to one decimal place.

\subsection{Tables definition}

Tables 1 through 4 refer to SP problems. Tables 1,2 and 3 refers to sets of 30 graphs with $|V| \in \lbrace 50,75,100,150 \rbrace$.  Table 4 refers to sets of 30 graphs with $|V| \in \lbrace 50,75,100\rbrace$  generated  independently of tables 1,2 and 3. There is one exception: in table 2 the line with $(|V|,K,\delta) = (100,20,0.25)$ refers to 5 problems. We are considering the interval case (tables 1 and 2), the line segment case (table 3) and the budgeted case (table 4) according the parameters presented. Partitions generation procedures are identified with $r,p$ or $d$.\\  

Tables 5 trough 7 refer to $(l,p)M$ problems with $30$ cases for each $l \in \lbrace 50,60,70,80,90 \rbrace$. In each case we are looking for $p = l/10$ medians. There are some exceptions: in table 5 lines with $( (l,p),\delta) = ((80,8),0.25)$ and $( (l,p),\delta) = ((90,9),0.25)$ refers to 10 problems, in table 6 line with $(l,p)=(90,9)$ refers to 10 problems and in table 7 lines in the right side with $((l,p),K) = ((80,8),15)$ refers to 15 cases. We are considering the interval case (table 5), the line segment case (table 6) and the budgeted case (table 7) according the parameters presented. Partitions generation procedures are identified with $lo$ or $g$.

\subsection{Performance of the algorithm}

The experimental results are far from being exhaustive. There are many factors involved: the dimensions of the problems ($|V|$ and $(l,p)$), the dimension of the partition ($K$), the definitions used for $P_k$ ($r,p,d$ for SP problems and $lo,g$ for $(l,p)$M problems), the definitions used for $\Omega$ (interval,line segment and budgeted) and the parameters that we need to define them ($\delta,\underline{ \alpha},\overline{\alpha},\beta_1,\beta_2)$.\\

We try to show that the multiparametric analysis to find an $\epsilon,\Omega$-mprs is computationally possible for moderate $|V|,l,p,K$ values and for this we define a variety of reasonable situations for the data used to define the cases. Remember that the algorithm will run offline and the decision maker will choose a solution from the generated set when a new scenario appears. Thus, a tolerable execution time may be hours and not minutes as usual. Also, the trade off analysis between the number of generated solutions and the relative error used defined with $\epsilon$ depends on the decision maker.\\

\subsubsection{General comments}

\begin{itemize}
	
\item As we can expect the computational effort and the number of generated solutions  increase if either $(|V|,K)$ or $((l,p),K)$ increase.

\item The known upper bound for the relative error at the first (or second) iteration may be very high which justify the use of more solutions.

\item Execution time seems to be tolerable (from seconds up to 6 hours). The number of generated solutions were low in some cases and tolerable in general although the trade off analysis between the number of solutions and the relative error depends on decision maker.

\end{itemize} 

\subsubsection{Remarks}

\begin{itemize}

\item We can see in tables 1,2 and 5 for the interval case that the computational effort and the number of generated solutions   increase if $\delta$ decreases. We can see in tables 2 and 5  that the known upper bound  for the relative error on the first iteration increases if $\delta$ decreases which is consistent with that behavior. 
In Appendix B we present an easy problem with $n+1$ variables called $Toy_n(c)$ to illustrate that behavior. For $Toy_n(c)$  we have that with low uncertainty (see Appendix B) we need $n$ solutions to find a $0,\Omega$-mprs and with high uncertainty (see Appendix B) one solution is enough to find a $0,\Omega$-mprs. That only suggests that finding $\epsilon,\Omega$-mprs may be easy with high uncertainty and may be hard with low uncertainty. We are not comparing the goodness of both situations. Doing so would surely lead to a redefinition of the objective: with low uncertainty it is most likely not necessary to be so picky about the value of $\epsilon$ and although we might be relatively far from robust solutions, in practice we might have reasonable solutions.

\item In table 2 we show that the decision maker can observe the increasing of the number of necessary solutions ($\overline{s_{d,q}}$) and the time ($\overline{t_{d,q}}$) to find an $(q/100) \times  v(R_{\mathcal{L}}),\Omega$-mprs while it is more demanding for stopping  the algorithm $(q \in \lbrace 5,3,2,1 \rbrace)$.

\item Note that for the previously defined budgeted case for $(l,p)M$ problems, if $\delta$ increases the volume of $\Omega$ increases with the interval case as its limit, hence the generated solutions must increase. However, the computational effort can increase and then decrease and, as expected, there are some cases with that behavior.

\end{itemize}

\section{Multiparametric analysis for a variant for the uncertainty set}

Let us suppose that the uncertainty set is an usual variant of $\Lambda(\Gamma)$ defined as follows:
\[\Lambda^+(\Gamma) = \lbrace c \in \mathbb{R}^n_+: c_j = \underline{c}_j + \lambda_j d_j\;\forall j \in [n],\;\underset{j \in P_k}{\sum} \lambda_j \leq \Gamma_k\;\forall k \in [K],\;\lambda \in [0,1]^n\rbrace\]

Next we present how to apply the multiparametric analysis for this variant. Problems $W(x,\Gamma)$ and $R(\Gamma)$ are now defined as follows: 
\begin{flalign*}  
&\underset{c \in \Lambda^+(\Gamma)}  {\max}\;\;c^t x&W(x,\Gamma)&\\
&\underset{x \in X} {\min}\;\;v(W(x,\Gamma)) = \underset{x \in X} {\min}\;\left ( \underset{c\in \Lambda^+(\Gamma)}{\max}\;c^t x \right )&R(\Gamma)\
\end{flalign*}

As before, dualizing the inner problem in $R(\Gamma)$ lead us to a 0-1-MILP problem in $(\pi,\rho,x)$ defined as follows (\cite{GoeLen2021}): 
\begin{flalign*}
\min&\;\;\Gamma^t \pi + \underset{j \in [n]}{\sum} \rho_j + \underline{c}^t  x&R(\Gamma)&\\
s.t.&\;\;\pi_k + \rho_j - d_j x_j \geq 0&\forall k \in [k]\;\forall j \in P_k&\\
&\;\;\pi \in \mathbb{R}^K_+,\;\rho \in \mathbb{R}^n_+,\;x \in X&\
\end{flalign*}

Unfortunately the parameters in $\Gamma$ are affecting a vector of continuous variables ($\pi$) and the multiparametric analysis by using {\bf A-$Q$} (see Appendix A) is not possible with the formulation presented. Next we present a new formulation.\\

Let $(\pi,\rho,x)$ be an optimal solution for $R(\Gamma)$ then because of the minimization criterium we have  $\rho_j = \max \lbrace 0,d_j x_j - \pi_k \rbrace\;\forall j \in P_k\;\forall k \in [K]$. Hence $R(\Gamma)$ may be rewritten as a problem in $(\pi,x)$ as follows:
\begin{flalign*}
\min&\;\;\underset{k \in [k]}{\sum} \left (\Gamma_k \pi_k + \underset{j \in P_k}{\sum} \max \lbrace 0,d_j x_j - \pi_k \rbrace \right ) + \underline{c}^t  x&R(\Gamma)&\\
s.t.&\;\;\pi \in \mathbb{R}^K_+,\;x \in X&\
\end{flalign*}

Let $(\pi^*,x)$ be an optimal solution for $R(\Gamma)$ then $\pi^*_k$ is an optimal solution for the following problem in $\pi_k$:
\begin{flalign*}
\min&\;\;\Gamma_k \pi_k + \underset{j \in P_k}{\sum} \max \lbrace 0,d_j x_j - \pi_k \rbrace &R(\Gamma_k,x)&\\
s.t.&\;\;\pi_k \in \mathbb{R}_+&\
\end{flalign*}

If $d_j x_j = 0\;\forall j \in [n]$ then $0$ is an optimal solution for $R(\Gamma_k,x)$. Otherwise the objective function in $R(\Gamma_k,x)$ is a convex, continuous and piecewise linear function  on $[0,\infty)$ (see Appendix C) and then we have that there exists an optimal solution with $\pi_k \in \lbrace 0 \rbrace \cup \lbrace d_j x_j:j \in P_k \rbrace$. Therefore $R(\Gamma)$ may be rewritten as a problem in $(\rho,x,\alpha)$ as follows:
\begin{flalign*}
\min&\;\;\underset{k \in [K]}{\sum} \left (\Gamma_k \underset{s \in P_k}{\sum} d_s \alpha_s  x_s \right )  + \underset{j \in [n]}{\sum} \rho_j + \underline{c}^t  x&R(\Gamma)&\\
s.t.&\;\;\underset{s \in P_k}{\sum} d_s \alpha_s x_s + \rho_j - d_j x_j \geq 0&\forall k \in [K]\;\forall j \in P_k&\\
&\;\;\underset{s \in P_k}{\sum} \alpha_s \leq 1&\forall k \in [K]&\\
&\;\;\alpha \in {\lbrace 0,1 \rbrace}^n,\;\rho \in \mathbb{R}^n_+,\;x \in X&\
\end{flalign*}

Let $(\rho,x,\alpha)$ be an optimal solution for $R(\Gamma)$. Let $\pi_k = \underset{s \in P_k}{\sum} d_s \alpha_s x_s\;\forall k \in [K]$ then $\pi_k \in \lbrace 0 \rbrace \cup \lbrace d_j x_j:j \in P_k \rbrace$ and the reformulation of $R(\Gamma)$ is valid.\\

$R(\Gamma)$ may be rewritten as a 0-1-MILP problem in $(\rho,x,\alpha,w)$ as follows:
\begin{flalign*}
\min&\;\;\underset{k \in [K]}{\sum} \left (\Gamma_k \underset{s \in P_k}{\sum} d_s w_s\right )  + \underset{j \in [n]}{\sum} \rho_j + \underline{c}^t  x&R(\Gamma)&\\
s.t.&\;\;\underset{s \in P_k}{\sum} d_s w_s + \rho_j - d_j x_j \geq 0&\forall k \in [K]\;\forall j \in P_k&\\
&\;\;\underset{s \in P_k}{\sum} \alpha_s \leq 1&\forall k \in [K]&\\
&\;\;w_j - \alpha_j -x_j \geq -1&\forall j \in [n]&\\
&\;\;w \in {\lbrace 0,1 \rbrace}^n,\;\alpha \in {\lbrace 0,1 \rbrace}^n,\;\rho \in \mathbb{R}^n_+,\;x \in X&\
\end{flalign*}

Let $(\rho,x,\alpha,w)$ be an optimal solution. If $x_j = 0$ then $w_j \geq \alpha_j - 1$ and because of the minimization  criterium  we have $w_j = 0 = \alpha_j x_j$. If $x_j = 1$ then $w_j \geq \alpha_j$ and because of the minimization criterium  we have $w_j = \alpha_j = \alpha_j x_j$ and the reformulation of $R(\Gamma)$ is valid. Now the parameters in $\Gamma$ are affecting only a vector of 0-1-variables ($w$), and then the mulltiparametric analysis may be performed by using {\bf A-$Q$} with problem $Q$ rewritten appropriately.

\section{Conclusions and further extensions}

In this paper we studied combinatorial problems with locally budgeted uncertainty parameterized with $\Gamma$. We presented an algorithm to find ${\lbrace x^i \rbrace}_1^r$ such that for any parameters vector $\Gamma$ there exists $i(\Gamma) \in [r]$ such that robustness of $x^{i(\Gamma)}$ is near optimal. As far as we know this is the first algorithm presented to do that task.\\  

The algorithm consists of applying a multiparametric algorithm to obtain a near optimal multiparametric solution relative to the objective function for a combinatorial problem defined to find a robust solution for the parameters vector fixed ($R(\Gamma)$). The case in which the matrix to define the nominal problems is unimodular is particulary considered.\\ 

Three cases for the parameters were considered for the computational experience: interval case, line segment case and budgeted case. Two problems were considered for the computational experience: the first one with the totally unimodular property,  the shortest path problem, and the second one without that property, the $p$-medians problem. The parameter to stop the algorithm was choice in such a manner that when the algorithm stops we have a near optimal solution in the relative sense.\\ 

The experience show that we may find ${\lbrace x^i \rbrace}_1^r$ with a tolerable computational effort for problems with moderate size  as we can see in tables. The decision maker can observe the increasing of the number of necessary solutions and computational effort according dimensions and uncertainty level.\\

For large problems the experience strongly suggests that we will need to solve the $Q$ problems (see {\bf A-$Q$}) more efficiently than using a standard branch and cut algorithm directly, either by using particular properties for the nominal problem or by using general properties for $Q$. The formulation for $R(\Gamma)$ that uses only the $\pi$ variables (see subsection 2.1) suggests that either a  relax (solve $\overline{Q}$ instead of $Q$ to obtain $\pi$ and a lower bound) and fix (solve $P(c(\pi))$ to obtain $(\rho,x)$ and an upper bound) with the framework of a branch and cut algorithm or a branch and cut algorithm  with the branching scheme guided with the $\pi$ variables are plausible options to be studied.\\

For low uncertainty, at least in the interval case, it may be necessary to redefine the approach if the number of solutions is not tolerable for the decision maker and perhaps a set with few solutions is enough if for any parameters vector under consideration some solution in the set  has a tolerable value. Note that to use a high $\epsilon$ value is not the unique option to be considered: {\it compromise solutions} to take to account the robustness relative error and the values at the same time may be a plausible option.\\

If the uncertainty set is parameterized either with $d$ instead of $\Gamma$ or $d$ and $\Gamma$ instead of $\Gamma$ the approach remains valid with a straightforward redefinition of $Q$, however we can expect a computational effort and a number of generated solutions not tolerable and then a redefinition of the objective may be necessary.\\ 

Our approach remains valid for all uncertainty sets such that there exists a 0-1-MILP formulation for the problem to find a robust solution with the fixed parameters affecting only 0-1-variables. That is the case for the variant of the uncertainty set presented.\\ 

\section*{Declaracion of interest} 

The author declares that he has no known competing financial interests or personal relationships that could have appeared to influence the work reported in this paper.

\section*{Acknowledgments}

Research supported by Universidad Central de Venezuela 

\bibliography{mybibfileEUROJCO2023}

\section*{Appendix A: multiparametric algorithm relative to the objective function for 0-1-ILP problems}

We present a summary of the results that may be seen in \cite{Cre2000},\cite{Cre2020}.\\

Let $X \subseteq {\lbrace 0,1 \rbrace}^n$ with $X \neq \emptyset$. Let us suppose that $X$ is 0-1-Mixed Integer Linear Programming (0-1-MILP) representable and let $P(c)$ be a combinatorial problem in $x$ parameterized in $c \in \mathbb{R}^n_+$, defined as follows:
\begin{flalign*}
\underset{x \in X}{\min}&\;\;c^t x&P(c)&\
\end{flalign*}

Let $\Omega \subseteq \mathbb{R}^n_+$, let $\mathcal{U} \in \mathbb{R}^n_+$ and let us suppose that $c \leq \mathcal{U}$ for all $c \in \Omega$.       
Let $\epsilon \geq 0$ and let ${\lbrace x^i \rbrace}_1^r \subseteq X$ such that:                           
$\underset{i \in [r]}{\min}\;c^t x^i-v(P(c)) \leq \epsilon\;\forall c \in \Omega$ then we say that ${\lbrace x^i \rbrace}_1^r$ is a $\epsilon,\Omega$-optimal multiparametric solution ($\epsilon,\Omega$-omps).\\

Let $Q({\lbrace x^i \rbrace}_1^r)$ be  a problem in $(c,x)$ defined as follows:
\begin{flalign*}
\max&\;\;\left (\underset{i \in [r]}{\min}\;c^t x^i-c^t x \right )&Q({\lbrace x^i \rbrace}_1^r)&\\
s.t.&\;\;c \in \Omega,\;x \in X&\
\end{flalign*}

We know that:

\begin{enumerate}
	\item If $v(Q({\lbrace x^i \rbrace}_1^r)) \leq \epsilon$ then ${\lbrace x^i \rbrace}_1^r$ is an $\epsilon,\Omega$-omps
	\item If $(c,x)$ is an optimal solution for $Q({\lbrace x^i \rbrace}_1^r)$ then $x$ is an optimal solution for $P(c)$
\end{enumerate}

{\bf Multiparametric algorithm ({\bf A-$Q$})}\\

Let $\epsilon \geq 0$. Let $c \in \Omega$. Solve $P(c)$, let $x^1$ be an optimal solution and let $r=1$. 

\begin{enumerate}  
	\item Solve $Q({\lbrace x^i \rbrace}_1^r)$ and let $x$ be an optimal solution
	\item If $v(Q({\lbrace x^i \rbrace}_1^r)) \leq \epsilon$ STOP
	\item Let $x^{r+1} =x$, let $r = r+1$ and return to step 1
\end{enumerate}

{\bf A-$Q$} is finite and if ${\lbrace x^i \rbrace}_1^r$ is the output then ${\lbrace x^i \rbrace}_1^r$ is an $\epsilon,\Omega$-omps.\\

In order to solve $Q({\lbrace x^i \rbrace}_1^r)$ we may rewritten as  a 0-1-MILP problem in $(c,x,w,\sigma)$  as follows:
\begin{flalign*}
\max&\;\;\sigma - \underset{j \in [n]}{\sum} w_j&Q({\lbrace x^i) \rbrace}_1^r)&\\
s.t.&\;\;\sigma - c^t x^i \leq 0&\forall i \in [r]&\\
&\;\;w_j - c_j - \mathcal{U}_j x_j\geq -\mathcal{U}_j&\forall j \in [n]&\\ 
&c \in \Omega,\;x \in X,\;w \in \mathbb{R}^n_+,\;\sigma \in \mathbb{R}&\
\end{flalign*}

Let $\mathcal{L} \in \mathbb{R}^n_+$ with $\mathcal{L} \leq \mathcal{U}$ and let $\Omega = [\mathcal{L},\mathcal{U}] = \lbrace c \in \mathbb{R}^n_+: \mathcal{L} \leq c \leq \mathcal{U} \rbrace$. We know that:
\begin{enumerate} 
	\item Let $x \in X$ and let $c^+(x)_j = \mathcal{L}_j x_j + \mathcal{U}_j (1-x_j)$ for all $j \in [n]$ then $x$ is an optimal solution for $P(c)$ for some $c \in \Omega$ if and only if $x$ is an optimal solution for $P(c^+(x))$
	\item Let $Q^+({\lbrace x^i \rbrace}_1^r)$ be a problem in $x$ defined as follows:
	\begin{flalign*}
	\max&\;\;\left (\underset{i \in [r]}{\min}\;{c^+(x)}^t x^i-\mathcal{L}^t x \right )&Q^+({\lbrace x^i \rbrace}_1^r)&\\
	s.t.&\;\;x \in X&\
	\end{flalign*}
	We know that: 
	\begin{enumerate}
		\item $v(Q^+({\lbrace x^i \rbrace}_1^r) = v(Q({\lbrace x^i \rbrace}_1^r)$
		\item If $(c,x)$ is an optimal solution for $Q({\lbrace x^i \rbrace}_1^r)$ then $x$ is an optimal solution for $Q^+({\lbrace x^i \rbrace}_1^r)$
		\item If $x$ is an optimal solution for $Q^+({\lbrace x^i \rbrace}_1^r)$ then $(c^+(x),x)$ is an optimal solution for $Q({\lbrace x^i \rbrace}_1^r)$ 
	\end{enumerate}
\end{enumerate}

In order to solve $Q^+({\lbrace x^i \rbrace}_1^r)$ we may rewritten as q  0-1-MILP in $(x,\sigma)$ problem  as follows:
\begin{flalign*}
\max&\;\;\sigma - \mathcal{L}^t x&Q^+({\lbrace x^i) \rbrace}_1^r)&\\
s.t.&\;\;\sigma + {f^i}^t x \leq \mathcal{U}^t x^i&\forall i \in [r]&\\
&\;\;x \in X,\;\sigma \in \mathbb{R}&\
\end{flalign*}

with $f^i_j = (\mathcal{U}_j -\mathcal{L}_j) x^i_j\;\forall j \in [n]\;\forall i \in [r]$\\

Note that the approach may be easily extended to a 0-1-MILP problem with the parameters affecting only the 0-1 variables.

\section*{Appendix B: Toy example} 

Let $Toy_2(c)$ be a combinatorial problem in $x$ defined as follows:
\begin{flalign*}
\min&\;\;c_1 x_1 + c_2 x_2 + c_3 x_3&Toy_2(c)&\\
s.t.&\;\;x_1 + x_2 + x_3 = 1&\\
&\;\;x_j \in \lbrace 0,1 \rbrace&\forall j \in [3]&\
\end{flalign*}

Let $\underline{c}_1 = \underline{c}_2 = 10$ and let $\underline{c}_3 = 11.5$. Let $d_1 = d_2 = 2$ and let $d_3 = 0$. Let $K = 3$ and let $P_k = \lbrace k \rbrace$ for all $k \in [3]$.\\ 

The uncertainty set is defined as follows:
\[\Lambda(\Gamma) = \lbrace c \in \mathbb{R}^3: c = \underline{c} + \lambda,\;\underset{j \in P_k}{\sum} \lambda_j \leq \Gamma_k\;\forall k \in [3],\;\lambda_j \in [0,d_j]\;\forall j \in [n]   \rbrace\]

Let $\mathcal{L}^1_j= 2\;\forall j \in [n]$ and let $\mathcal{U}^1_j = 3\;\forall j \in [n]$. Let $\Omega^1 = [\mathcal{L}^1,\mathcal{U}^1]$. If $\Gamma \in \Omega^1$ then:\\

\begin{enumerate}
	\item  If $x^t = (1,0,0)$ then $v( W ( x,\Gamma) ) = 12$ and if $y^t = (0,1,0)$ then  $v(W(y,\Gamma)) = 12$\\
	\item  If $z^t = (0,0,1)$ then $v(W(z,\Gamma)) = 11.5$.
\end{enumerate}

Therefore $z$ is a robust solution for all $\Gamma \in \Omega^1$\\

Let $\mathcal{L}^2 = 0$ and let $\mathcal{U}^2_j = 1\;\forall j \in [n]$. Let $\Omega^2 = [\mathcal{L}^2,\mathcal{U}^2]$. If $\Gamma \in \Omega^2$ then:\\

\begin{enumerate}
	\item If $x^t = (0,0,1)$ then  $v(W(x,\Gamma)) = 11.5$ for all $\Gamma \in \Omega^2$\\
	
	\item  Let $y^t = (1,0,0)$ and $z^t = (0,1,0)$ then:
	
	\begin{enumerate}	
		\item  If $\Gamma_1 < \Gamma_2$ then $v(W(y,\Gamma)) = 10 + \Gamma_1 < 10 + \Gamma_2 = v(W(z,\Gamma))$ with $10 + \Gamma_1 < 11.5$.\\
		
		\item  If $\Gamma_2 < \Gamma_1$ then $v(W(z,\Gamma)) = 10 + \Gamma_2 < 10 + \Gamma_1 = v(W(y,\Gamma))$ with $10 + \Gamma_2 < 11.5$.\\
		
		\item  If $\Gamma_1 = \Gamma_2$ then $v(W(y,\Gamma)) = v(W(z,\Gamma)) = 10 + \Gamma_1 = 10 + \Gamma_2 < 11.5$
		
	\end{enumerate}	
\end{enumerate}

Therefore: $y$ is a robust solution if and only if $\Gamma_1 \leq \Gamma_2$ and $z$ is a robust solution if and only if $\Gamma_2 \leq \Gamma_1$ and we need two robust solutions to find a $0,\Omega^2$-mprs.\\

Next, we present a generalization. Let $Toy_n(c)$ be a combinatorial problem in $x$ defined as follows:
\begin{flalign*}
\min&\;\;c_1 x_1 +  \cdots + c_{n+1} x_{n+1}&Toy_n(c)&\\
s.t.&\;\;x_1 + \cdots + x_{n+1} = 1&\\
&\;\;x_j \in \lbrace 0,1 \rbrace&\forall j \in [n+1]&\
\end{flalign*}

Let $\underline{c}_1 = \cdots = \underline{c}_n = 10$ and let $\underline{c}_{n+1} = 11.5$. Let $d_1 = \cdots = d_n = 2$ and let $d_{n+1} = 0$. Let $K = n+1$ and let $P_k = \lbrace k \rbrace$ for all $k \in [n+1]$.\\

Let $\mathcal{L}^1_j= 2\;\forall j \in [n]$ and let $\mathcal{U}^1_j = 3\;\forall j \in [n]$. Let $\Omega^1 = [\mathcal{L}^1,\mathcal{U}^1]$. Let $\mathcal{L}^2 = 0$ and let $\mathcal{U}^2_j = 1\;\forall j \in [n]$. Let $\Omega^2 = [\mathcal{L}^2,\mathcal{U}^2]$.\\

Now the number of solutions to define a  $0,\Omega^2$-mprs  is $n$ and the number of solutions to define a  $0,\Omega^1$-mprs is 1. However, the uncertainty in $\Omega^1$ is larger than the uncertainty in $\Omega^2$ and the volume defined as $\underset{k \in [n+1]}\Pi (\mathcal{U}^s_k -\mathcal{L}^s_k)$ is the same for $s \in [2]$

\section*{Appendix C: Properties of an auxiliary function}

Let $v \in \mathbb{R}^q$ with $v > 0$, let $u \in \mathbb{N}^q$, let $\gamma \in \mathbb{R}_+$ and let $f(\pi) = \gamma \pi + \underset{j \in [q]}{\sum} u_j \max \lbrace 0,v_j - \pi \rbrace$ for all $\pi \in \mathbb{R}_+$ then $f$ is a continuous, convex and piecewise linear function. Next we present a proof.\\

Let $v_0 = 0$ and let us suppose without loss of generality that $v_0  < v_1 < \cdots < v_q$. Note that the $q+1$ pieces of $f$ are defined as follows:

\begin{itemize}
	\item $f(\pi) = f_k(\pi) = \gamma \pi + \overset{q}{\underset{j=k}{\sum}} u_j (v_j - \pi) = \overset{q}{\underset{j=k}{\sum}} u_j v_j + (\gamma - \overset{q}{\underset{j=k}{\sum}} u_j) \pi$ for $\pi \in [v_{k-1},v_k)$ for $k \in [q]$
	\item $f(\pi) = f_{q+1}(\pi) = \gamma \pi $ for $\pi \in [v_q,\infty)$
\end{itemize}

Note that $f_{k+1}(v_k) = \overset{q}{\underset{j=k+1}{\sum}} u_j v_j + (\gamma - \overset{q}{\underset{j=k+1}{\sum}} u_j) v_k = \overset{q}{\underset{j=k}{\sum}} u_j v_j + (\gamma - \overset{q}{\underset{j=k}{\sum}} u_j) v_k$ for $k \in [q-1]$ and note that $f_{q+1}(v_q) = \gamma v_q = u_q v_q + (\gamma- u_q) v_q$, therefore $f$ is a continuos piecewise function on $[0,\infty)$.\\

Note that successive slopes satisfy $(\gamma - \overset{q}{\underset{j=1}{\sum}} u_j) < (\gamma - \overset{q}{\underset{j=2}{\sum}} u_j) < \cdots < \gamma - u_q < \gamma$, therefore $f$ is convex.

\section*{Appendix D: Remark about the computational complexity}  

$R(\Gamma)$ may be rewritten as follows (\cite{GoeLen2021}):
\begin{flalign*}
\min&\;\;\Gamma^t \pi + v(P(c(\pi)))&\hat R(\Gamma)&\\
s.t.&\;\;\pi \in {\lbrace 0,1 \rbrace}^K&\
\end{flalign*}

Let $\hat X \subseteq X$ such that $x \in \hat X$ if and only there exists $\pi \in {\lbrace 0,1 \rbrace}^K$ such that $x$ is an optimal solution for $P(c(\pi))$. We have that $\hat X$ is a $0,\Omega$-mprs. In order to obtain $\hat X$ we need to solve $2^K$ nominal problems. Therefore if $P(c)$ may be solved in polynomial time and $K$ is constant we can obtain $\hat X$ in polynomial time.\\

\newpage

\begin{table} 
	\caption{SP. Interval case }
	\label{tab:4}
	\begin{tabular}{llllllllllllll}
		\hline\noalign{\smallskip}
		$|V|$ & $K$ & $\delta$ & $\overline{t_r}$&  $\hat t_r$ & $\underline{s_r}$ & $\overline{s_r}$ & $\hat s_r$ & $\overline{t_p}$&  $\hat t_p$ &$\underline{s_p}$ &$\overline{s_p}$   & $\hat s_p$  \\
		\noalign{\smallskip}\hline\noalign{\smallskip}
		50  & 5 &  0.0 &    2.2 &  3.3 & 8 & 15.7 &  22 & 2.3 &   4.5 & 4 & 7.1 & 12\\
			&   &  0.5 &   1.0 &  1.7 & 2 &  7.0 &  12 & 2.5 &   5.0 & 2 & 6.3 & 11\\
			&   &  1.0 &   0.7 &  1.1 & 1 &  4.4 &   7 & 0.8 &   1.7 & 1 & 3.3 &  8\\
			&   &  1.5 &  0.8 &  2.0 & 1 &  4.3 &   6 & 0.3 &   0.4 & 1 & 1.0 &  1\\
			& 10&  0.0 & 59.0 &109.2 & 45 & 108.6 &  204 & 84.6 &   211.9 &20 &38.9 & 80\\       
			&   &  0.5 &  4.7 & 10.8 & 2 & 18.8 &  33 & 14.3 &  30.0 & 6 &12.7 & 26 \\ 	
			&   &  1.0 &  1.8 &  3.8 & 2 &  7.8 &  13 & 0.4 &   1.0 & 1 & 1.4 &  3\\ 	
			&   &  1.5&  1.3 &  2.2 &  2 &  5.3 &    8 &  0.3 &   0.3 & 1 & 1.0 &  1\\
		\noalign{\smallskip}\hline
	\end{tabular}\\
\end{table}

\newpage

\begin{table} 
	\caption{SP. Interval case. *: refers to 5 cases.}
	\label{tab:4}
	\begin{tabular}{llllllllllllll}
		\hline\noalign{\smallskip}
	$|V|$ & $K$ & $\delta$ &$\overline{\epsilon_{\mathcal{L}}}$& $\underline{s_{d,1}}$  &   $\hat{s_{d,1}}$ &$\overline{s_{d,1}}$ &$\overline{s_{d,2}}$ & $\overline{s_{d,3}}$ & $\overline{s_{d,5}}$ &$\overline{t_{d,1}}$ & $\overline{t_{d,2}}$ & $\overline{t_{d,3}}$ & $\overline{t_{d,5}}$ \\
		\noalign{\smallskip}\hline\noalign{\smallskip}
		 75   & 10   &   0.25 & 40.9 & 23 & 73 & 43.7 & 39.6 & 36.5 & 30.1 &  117.8 &  97.7 &  81.3 & 51.9\\
		    &     &  0.50  & 34.0&12 & 44 & 23.2 & 21.3 & 19.5 & 16.4 &  53.3 &  44.8 &  37.9 & 24.9\\  
		    &   &  0.75 & 27.0 & 2 & 18 &  9.8 &  8.9 &  8.2 &  7.5 &  16.4 &  12.5 &  10.4 &  8.6\\
		    &   &     1.00 & 9.4&  1 &   6 &  2.3 &  2.2 &  2.2 &  2.2 &    2.1 &   1.9 &   1.9 &  1.9 \\
		    & 15  &  0.25 &  44.5 & 40 & 152 &96.8 & 89.2 & 82.8 & 69.7 &  645.3 & 522.6 & 437.2 & 260.7\\ 	
		    &   &  0.50  & 36.7&19 & 72 & 39.8 & 36.7 & 33.2 & 28.8 & 187.2 & 148.6 & 115.9 & 80.0\\
		    &   &  0.75 &  27.8 &3 & 25 & 11.4 & 10.7 &  9.7 &  8.9 &  34.2 &  26.1 &  20.7 & 17.7\\
		    &   &   1.00    & 5.8 & 1 &   4 &  1.6 &  1.6 &  1.6 &  1.6 &    2.3 &   2.3 &   2.3 & 2.3\\
100		 & 15 &  0.25 & 42.2&  77 & 224 &155.3 &141.8 &130.2 &105.2 & 3115.6 &2461.7 &1948.0 &1121.3 \\    
		    &  &  0.50  & 38.0&23 &78 &48.9 &44.2 &40.5 &35.5 &543.0 &429.0 &363.5 &256.2\\
		    &   &  0.75 & 28.0 & 8 &27 &16.3 &14.6 &13.4 &12.1 &131.1 &101.0 &80.8 &67.8\\
		    &   &  1.00  & 7.2  & 1 & 3  &1.8  &1.8  &1.8  &1.8  & 7.6   &7.4   &7.4  &7.4\\
		  & 20   &  0.25*   & 44.0 & 132 & 391 &259.2 &231.2 &207.4 &171.0 & 4754.4 &3869.5 &3180.9 &1914.0\\
			&    &  0.50 &39.1  &32 & 155 & 76.4 & 69.0 & 63.8 & 54.3 & 1428.9 & 1118.5 & 918.6 & 629.4 \\
			&   &  0.75 &  28.8&7 & 33& 19.2 & 16.7 & 15.5& 14.4& 222.8& 150.6& 124.1& 112.1\\
			&   &   1.00    & 2.4& 1 & 2 &1.2 &1.2 &1.2 &1.2 & 7.6 &7.6 &7.6 &7.6\\
		\noalign{\smallskip}\hline
	\end{tabular}\\ 
\end{table}

\newpage

\begin{table} 
	\caption{SP. Line segment case. If $|V| \in \lbrace 50,75 \rbrace$ then $*=p$. If $|V| \in \lbrace 100,150 \rbrace$ then $*=d$ }
	\label{tab:4}
	\begin{tabular}{llllllllllllll}
		\hline\noalign{\smallskip}
	$|V|$ & $K$ & $\overline{\epsilon_0}$ &$\overline{t_r}$& $\hat t_r$  & $\underline{s_r}$ & $\overline{s_r}$ & $\hat s_r$ &$\overline{\epsilon_0}$ &$\overline{t_*}$& $\hat t_*$  & $\underline{s_*}$ & $\overline{s_*}$ & $\hat s_*$\\
		\noalign{\smallskip}\hline\noalign{\smallskip} 
		50  & 5 & 51.3&0.4 & 0.8 & 2 & 2.9 & 4 &31.2 &0.6 & 0.8 & 1 & 2.3 & 4\\
		    & 10& 119.2&0.6 & 0.9 & 1 & 3.2 & 5 &96.6 &1.2 & 2.0 &1 &2.3 &4\\
		75  & 10 & 123.0 &1.4 & 2.4 & 3 & 3.9 & 6 &98.4 &3.1 & 5.3 & 1 & 3.1 & 5\\
		    & 15 &192.4 &2.8 & 5.3 & 3 & 4.1 & 5 & 164.1 &6.8 & 17.3 & 1 & 2.9 & 5\\
		    			 100& 15& 190.1 &5.5 &10.7 &3 &4.1 &5 & 163.1&16.2& 27.3 &1 &3.2 &5\\
		    & 20 & 259.9&8.3 &14.5 &3 &4.3 &6 &226.3 &24.9& 39.0 &1 &3.4 &6\\
		 150 & 20 & 261.3 &21.1 & 33.8 & 3 & 4.1 & 5 & 228.1&99.8 &178.2 &2 &4.2 &6 \\    
		     & 25 & 330.2&29.4 & 53.1 & 3 & 4.2 & 5& 292.7&144.6 & 236.5 & 1 & 4.0 & 7\\         
		\noalign{\smallskip}\hline
	\end{tabular}\\
\end{table}	

\newpage

\begin{table} 
	\caption{SP. Budgeted case. $\beta_1=0.5,\;\beta_2=1.0$}
	\label{tab:4}
	\begin{tabular}{llllllllllllllll}
		\hline\noalign{\smallskip}
		$|V|$ & $K$ & $\delta$ & $\overline{\epsilon_{\underline{\Gamma}}}$ &$\overline{t_r}$& $\hat t_r$ & $\underline{s_r}$ & $\overline{s_r}$ & $\hat s_r$ & $\overline{\epsilon_{\underline{\Gamma}}}$ & $\overline{t_d}$& $\hat t_d$ &$\underline{s_d}$ &$\overline{s_d}$ & $\hat s_d$  \\
		\noalign{\smallskip}\hline\noalign{\smallskip}
		50  & 10& 2 & 7.9&0.8 & 1.5 & 1 & 6.0 & 11 & 8.8&2.4 & 7.3 & 2 & 8.0 & 18\\
		    &   & 3 & &1.0 &2.2 &1 &7.3 &15 & &3.8 &11.9 & 3 &10.4 &26\\
		    &   & 5 & &1.2 &3.3 &1 &8.1 &19 & &4.8 &15.2 &3 &11.1 &26\\
		75  & 10& 2 & 8.0 &2.0 &3.8 &3 &8.1 &14 & 9.3 &8.6 & 17.8 & 5 & 11.1 &22\\
		    &   & 3 & &2.8 & 6.1 & 3 & 10.3 & 22 & &13.8 &32.8 &6 &15.0 &30\\
		    &   & 5 & &3.4 & 8.9 & 3 & 12.1 & 28 & &19.4 & 47.2 & 9 & 17.0 & 33\\
		    & 15& 3 & 8.1 &2.8 &6.7 &3 &10.3 &22 & 12.3 &14.1 &32.5 &6 &15.0 &30\\
		    &   & 5 & &9.0 & 21.6 &3 &17.9 &30 & &72.1 &189.2 &15 &32.5 &60\\
		    &   & 8 & &11.3 &31.4 &3 &19.4 &34 & &92.9 &255.7 &15 &34.9& 69 \\
		100 & 15& 3 & 7.6&11.7 &26.1 &2 &10.4 &17 & 14.3 &104.5 &186.1 &19 &27.1 &38\\
		    &   & 5& &14.5 &39.9 &2 &14.2 &34& &190.6 &354.4 &23 &36.9 &57\\
		    &   & 8& &16.4 &47.2 &2 &15.5 &39& &229.8 &464.7 &22 &39.3  &56\\
		    & 20& 4 & 9.3&26.8 & 63.8 &3 &18.3 &36 &17.4 &417.4 & 834.1 & 24 & 47.3 & 78\\
		    &   & 7 & &44.9 & 135.0 &3 &25.2 &62 & &802.3 & 1596.4 &35 &67.6 &115\\
		    &   & 10 & &49.1 & 143.3 & 3 & 26.3 & 64 & &856.5 & 1750.0 & 36 & 69.8& 125\\           
		\noalign{\smallskip}\hline
	\end{tabular}\\
\end{table}

\newpage

\begin{table} 
	\caption{$(l,p)$M. Interval case. $K = 10$. *: refers to 10 cases   }
	\label{tab:4}
	\begin{tabular}{lllllllllllllll}
		\hline\noalign{\smallskip}
		$l,p$ & $\delta$ & $\overline{\epsilon_{\mathcal{L}}}$&$\overline{t_{lo}}$& $\hat t_{lo}$ & $\underline{s_{lo}}$ & $\overline{s_{lo}}$ & $\hat s_{lo}$ & $\overline{\epsilon_{\mathcal{L}}}$ &$\overline{t_g}$& $\hat {t_g}$ &$\underline{s_g}$ &$\overline{s_g}$ & $\hat {s_g}$  \\
		\noalign{\smallskip}\hline\noalign{\smallskip}
		50,5  &     0.25 & 51.2 &118.2 &290.7 &2 &9.8 &19 & 49.4&137.0 &505.1 &3 &11.0& 24 \\
		60,6  &          & 46.3&253.7 &1103.6 &1 &10.0 &20 & 41.5&314.2 &1007.1 &2 &10.9 &23\\  
		70,7  &          &46.5 &788.1 &3913.8 &4 &11.5 &27 & 49.4&964.3 &5720.6 &4 &12.0 &21 \\
		80,8* &     &53.4 &1175.6 &2363.9 &2 &10.2 &14 & 50.4&2231.3 &5986.8 &2 &11.0 &20 \\
		90,9*  & & 50.1 &5250.0 &15760.3 &6 &12.2 &23 & 53.3 & 6291.7 &20441.4 &2 &11.1 &20\\
		50,5  &     0.50 & 20.0 &33.1 &88.6 &1 &3.4 &9 & 25.4 &42.3 &214.2 &1 &4.0 &13 \\
		60,6  &          & 19.9 &69.9 &271.3 &1 &3.5 &8 & 18.9 &71.3 &227.5 &1 &3.4 &9 \\
		70,7  &          & 20.8&174.2 &787.0 &1 &3.6 &9 & 18.5&160.1 &879.8 &1 &3.0 &6\\ 
		80,8  &          &20.3 &328.9 &1081.3 &1 &2.7 &6 &16.8&388.0 &1220.2 &1 &2.9 &6 \\  
		90,9  &          &18.2 &915.6 &5500.3 &1 &2.6 &6 & 16.4&957.8 &5773.6 &1 &2.5 &9\\
		\noalign{\smallskip}\hline
	\end{tabular}\\
\end{table}

\newpage

\begin{table} 
	\caption{$(l,p)$M. Line segment case. $K = 10$. *: refers to 10 cases}
	\label{tab:4}
	\begin{tabular}{lllllll}
		\hline\noalign{\smallskip}
		$l,p$  & $\overline{{\epsilon_0^{(2)}}}$ & $\overline{t_g}$& $\hat {t_g}$ &$\underline{s_g}$ &$\overline{s_g}$ & $\hat {s_g}$  \\
		\noalign{\smallskip}\hline\noalign{\smallskip}
		50,5 &  31.4 & 42.7 &170.0 &1 &2.7 &4 \\
		60,6 &  30.5 &193.7 &772.6 &1 &2.9 &5\\  
		70,7 &  29.3 & 655.4 &1563.4 &1 &3.4 &5  \\
		80,8 &  33.7 &1476.1 &2784.4 &2 &3.8 &6  \\
		90,9* & 27.9 &3098.0 &5375.5 &1 &3.5 &6\\
		\noalign{\smallskip}\hline
	\end{tabular}\\
\end{table}

\newpage

\begin{table} 
	\caption{$(l,p)$M. Budgeted case. $\beta_1=0.25,\;\beta_2=1.0$. *:refers to 15 cases}
	\label{tab:4}
	\begin{tabular}{llllllllllllllll}
		\hline\noalign{\smallskip}
		$l,p$ & $K$ & $\delta$ & $\overline{\epsilon_{\underline{\Gamma}}}$ &$\overline{t_{lo}}$& $\hat t_{lo}$ & $\underline{s_{lo}}$ & $\overline{s_{lo}}$ & $\hat s_{lo}$ & $\overline{\epsilon_{\underline{\Gamma}}}$ & $\overline{t_g}$& $\hat t_g$ &$\underline{s_g}$ &$\overline{s_g}$ & $\hat s_g$  \\
		\noalign{\smallskip}\hline\noalign{\smallskip}
	70,7 & 10 & 2 & 8.2& 346.2 &2118.8 &1 &5.6 &11& 7.8 & 243.3 &885.3 &1 &4.8 &11 \\
	    &    &  3&8.5 & 511.8 &5168.7 &1 &6.4 &14& 8.4 &333.2 &1515.4 &1 &5.6 & 15\\
		&    &  5& 8.5&551.4 &5644.5 &1 &6.6 &16 &8.4 & 359.3 &1974.7 &1 &5.7 &16\\
		& 15 &  3& 9.4& 972.0 &8832.9 &2 &7.5 &19& 8.8& 546.8 &2023.7 &1 &7.2 &18\\
		&    &  5& 9.4 &1010.0 &7619.9 &2 &8.1 &22 & 8.8& 540.7 &1858.7 &1 &7.3 &18\\
		&    &  8& 9.4 &1036.3 &7461.1 &2 &8.1 &24 &8.8 & 484.8&1792.8 &1 &7.3 &18\\   
	80,8 & 10 & 2&9.1 & 1265.2& 4732.0& 1& 6.0& 13& 9.4 & 1453.3 &5592.1 &1 &6.5 &12 \\ 
		     &    &  3& 9.5&1397.3 &5315.4 &1 &6.3 &14 & 9.1& 1795.7 &6514.8 &1 & 7.1 & 12\\
		     &    &  5& 9.5&1323.6 &5640.1 &1 &6.5 &15  & 9.1& 1964.6 &8634.2 &1 &7.1 &13\\ 
		     & 15 &  3& 8.0 & 1498.8 &7760.3 &1 &5.9 &14 &*8.2  & 2090.4 &13435.4 &2 &7.3 &21 \\
		     &    &  5 & 8.0 &1585.4 &9399.7 &1 &6.2 &15&*8.2& 2806.4 &17669.0 &2 &7.7 &22 \\
		     &    &  8 & 8.0 & 1618.4 &9318.8 &1 &6.1 &15 &*8.2 & 2708.1 &13585.1 &2 &7.7 &22 \\
                  
		\noalign{\smallskip}\hline
	\end{tabular}\\
\end{table}

\end{document}